\documentclass[11pt]{article}

\usepackage{amssymb}
\usepackage{amsmath}
\usepackage{amsthm}
\usepackage{mathrsfs}
\usepackage{graphicx}
\usepackage[utf8]{inputenc}
\usepackage{multirow}
\usepackage{floatrow}
\usepackage{mathrsfs}

\newtheorem{thm}{Theorem}[section]
\newtheorem{prop}[thm]{Proposition}
\newtheorem{lem}[thm]{Lemma}
\newtheorem{cor}[thm]{Corollary}
\newtheorem{rem}[thm]{Remark}

\newtheorem{notn}[thm]{Notation}

\newenvironment{Pf}{\noindent{\bf Proof. }}{\hfill $\blacksquare$ \\}

\allowdisplaybreaks

\def\lf{\lfloor}
\def\rf{\rfloor}

\def\ov{\overline}

\def\la{\langle}
\def\ra{\rangle}

\def\a{\alpha}

\def\e{\textbf e}

\def\g{\textbf g}
\def\k{\kappa}
\def\l{\lambda}
\def\m{\textbf m}
\def\r{\rho}
\def\s{\sigma}

\def\Ap{\text{Ap}}
\def\F{\textbf F}
\def\G{\text G}

\def\PF{\text{PF}}

\def\Z{\mathbb Z}

\def\ds{\displaystyle}
\def\sm{\setminus}
\def\pr{\prime}

\begin{document}

\title{\bf On the Frobenius Problem for Some Generalized Fibonacci Subsequences - I}
\author{
{\bf Santak Panda}\thanks{Department of Mathematics, Ohio State University, Columbus, OH 43210-1174, USA. \newline {\tt e-mail:panda.58@osu.edu}} \qquad
{\bf Kartikeya Rai}\thanks{Institut de Math\'{e}matique d'Orsay, Universit\'{e} Paris-Saclay, F-91405 Orsay Cedex, France. \newline {\tt e-mail:rai.kartikeya@etu-upsaclay.fr}} \qquad
{\bf Amitabha Tripathi}\thanks{Department of Mathematics, Indian Institute of Technology, Hauz Khas, New Delhi -- 110016, India. \newline {\tt e-mail:atripath@maths.iitd.ac.in}}\:\:\thanks{\it Corresponding author}
}
\date{}
\maketitle

\begin{abstract}
\noindent For a set $A$ of positive integers with $\gcd(A)=1$, let $\la A \ra$ denote the set of all finite linear combinations of elements of $A$ over the non-negative integers. The it is well known that only finitely many positive integers do not belong to $\la A \ra$. The Frobenius number and the genus associated with the set $A$ is the largest number and the cardinality of the set of integers non-representable by $A$. By a generalized Fibonacci sequence $\{V_n\}_{n \ge 1}$ we mean any sequence of positive integers satisfying the recurrence $V_n=V_{n-1}+V_{n-2}$ for $n \ge 3$. We study the problem of determining  the Frobenius number and genus for sets $A=\{V_n, V_{n+d}, V_{n+2d}, \ldots \}$ for arbitrary $n$, where $d$ odd or $d=2$. 
\end{abstract}

\noindent {\bf Keywords.} Embedding dimension, Ap\'{e}ry set, Frobenius number, Genus

\noindent {\bf 2020 MSC.} 11D07, 20M14, 20M30

\section{Introduction} \label{intro}
\vskip 10pt

For a given subset $A$ of positive integers with $\gcd(A)=1$, we write
\[ S = \la A \ra = \big\{ a_1x_1+\cdots+a_kx_k: a_i \in A, x_i \in {\Z}_{\ge 0} \big\}. \]
We say that $A$ is a set of generators for the set $S$. Further, $A$ is a minimal set of generators for $S$ if no proper subset of $A$ generates $S$. Let $A=\{a_1, \ldots, a_n\}$ be a set of generators of $S$ arranged in increasing order. Then the following are equivalent: 
\begin{itemize}
\item $A$ is a minimal set of generators for $S$; \\[-16pt]
\item $a_{k+1} \notin \la a_1, \ldots, a_k \ra$ for $k \in \{1,\ldots,n-1\}$; \\[-16pt] 
\item $A=S^{\star} \sm \big(S^{\star}+S^{\star}\big)$, where $S^{\star}=S \sm \{0\}$. \\[-16pt]
\end{itemize}
The embedding dimension $\e(S)$ of $S$ is the size of the minimal set of generators. 

For any set of positive integers $A$ with $\gcd(A)=1$, the set ${\Z}_{\ge 0} \sm S$ is necessarily finite; we denote this by $\G(S)$. The cardinality of $\G(S)$ is the genus of $S$ and is denoted by $\g(S)$. The largest element in $\G(S)$ is the Frobenius number of $S$ and is denoted by $\F(S)$.

The Ap\'{e}ry set of $S$ corresponding to any fixed $a \in S$, denoted by $\Ap(S,a)$, consists of those $n \in S$ for which $n-a \notin S$. Thus, $\text{Ap}(S,a)$ is the set of minimum integers in $S \cap \textbf {C}$ as $\textbf{C}$ runs through the complete set of residue classes modulo $a$. 

The integers $\F(S)$ and $\g(S)$, and the set $\PF(S)$, can be computed from the Ap\'{e}ry set $\Ap(S,a)$ of $S$ corresponding to any $a \in S$ via the following proposition. 

\begin{prop} {\bf (\cite{BS62,Sel77,Tri06})} \label{prelims}
Let $S$ be a numerical semigroup, let $a \in S$, and let $\Ap(S,a)$ be the Ap\'{e}ry set of $S$ corresponding to $a$. Then 
\begin{itemize}
\item[{\rm (i)}]
\[ \F(S) = \max \Big( \Ap(S,a) \Big) - a; \]
\item[{\rm (ii)}]
\[ \g(S) = \frac{1}{a} \left( \sum_{n \in \Ap(S,a)} n \right) - \frac{a-1}{2}; \]
\item[{\rm (iii)}]
\[ \PF(S) = \Big\{ n-a: n \in \Ap(S,a), n+\m(x)>\m(x+n), x=1,\ldots,a-1 \Big\},\]
where $\m(x) \in \Ap(S,a)$ and $\m(x) \equiv x\pmod{a}$.  
%\item[{\rm (iv)}]
%If $A=\{a_1, \ldots, a_n\}$ is a set of generators of $S$ arranged in increasing order, then $A$ is a minimal set of generators if and only if $a_{k+1} \notin \la a_1, \ldots, a_k \ra$ for $1 \le k \le n-1$. 
\end{itemize}  
\end{prop}
\vskip 5pt

The case where $\e(S)=2$ is well known and easy to establish. If $S=\la a,b \ra$, then it is easy to see that $\Ap(S,a)=\{bx: 0 \le x \le a-1\}$, and consequently 
\begin{equation} \label{e(S)=2}
\F(S) = ab-a-b, \quad \g(S) = \tfrac{1}{2} (a-1)(b-1), \quad \PF(S) = \{ab-a-b\}
\end{equation} 
by Proposition \ref{prelims}. 
\vskip 5pt

The Frobenius Problem is the problem of determining the Frobenius number and the genus of a given numerical semigroup, and was first studied by Sylvester, and later by Frobenius; see \cite{Ram05} for a survey of the problem. Connections with Algebraic Geometry revived interest in Numerical Semigroups around the middle of the twentieth century; we refer to \cite{RG-S09} as a basic textbook on the subject. Curtis \cite{Cur90} proved that there exists no closed form expression for the Frobenius number of a numerical semigroup $S$ with $\e(S)>2$. As a consequence, a lot of research has focussed on the Frobenius number of semigroups whose generators are of a particular form. There are three particular instances of such results that are perhaps the closest to our work, and hence bear mentioning. Mar\'{i}n et. al.\cite{MRR07} determined the Frobenius number and genus of numerical semigroups of the form $\la F_i, F_{i+2}, F_{i+k} \ra$, where $i, k \ge 3$. These are called Fibonacci semigroups by the authors. Matthews \cite{Mat09} considers semigroups of the form $\la a, a+b, aF_{k-1}+bF_k \ra$ where $a>F_k$ and $\gcd(a,b)=1$. Taking $a=F_i$ and $b=F_{i+1}$, one gets the semigroup $\la F_i, F_{i+2}, F_{i+k} \ra$, considered in \cite{MRR07}. Thus, such semigroups were termed generalized Fibonacci semigroups by Matthews, who determined the Frobenius number of a generalized Fibonacci semigroup, thereby generalizing the result in \cite{MRR07} for Frobenius number. Batra et. al. \cite{BKT15} determined the Frobenius number and genus of numerical semigroups of the form $\la a, a+b, 2a+3b, \ldots, F_{2k-1}a+F_{2k}b \ra$ and $\la a, a+3b, 4a+7b, \ldots, L_{2k-1}a+L_{2k}b \ra$ where $\gcd(a,b)=1$. 

By a generalized Fibonacci sequence we mean any sequence $\{V_n\}$ of positive integers which satisfies the recurrence $V_n=V_{n-1}+V_{n-2}$ for each $n \ge 3$. This paper studies the numerical semigroup generated by some subsequences of a generalized Fibonacci sequence $\{V_n \}$. More specifically, we study the numerical semigroup $S$ generated by $\la V_n, V_{n+d}, V_{n+2d}, \ra$, where $V_1=a$ and $V_2=b$. Our main results are: 
\begin{itemize}
\item[{\rm (i)}] $S=\la V_n, V_{n+d}, V_{n+2d}, \ldots \ra$ is a numerical semigroup if and only if $\gcd(a,b)=1$ and $\gcd(V_n,F_d)=1$; see Theorem \ref{existence}.
\item[{\rm (ii)}] If $d$ is odd, then $\e(S)=2$. Consequently, the computation of $\F(S)$ and $\g(S)$ is straightforward; see Theorem \ref{e(S)=2} and Theorem \ref{odd_case}. 
\item[{\rm (iii)}] If $d=2$, then $\e(S)=\k$ where $\k$ satisfies $F_{2(\k-1)} \le V_n-1 < F_{2\k}$; see Theorem \ref{embed_dim}.  
\item[{\rm (iv)}] If $d=2$, then $\Ap(S,V_n)=\{ s(x): 1 \le x \le V_n-1 \} \cup \{0\}$, where $s(x)$ is obtained by applying the Greedy Algorithm to $x$ with respect to the sequence $F_2, F_4, F_6, \ldots $; see Theorem \ref{Apery_set}.  
\item[{\rm (v)}] If $d=2$, then $\g(S)=s(V_n-1)-V_n$ and $\g(S)=\frac{1}{V_n} \left( \sum_{x=1}^{V_n-1} s(x) \right) - \frac{V_n-1}{2}$; see Theorem \ref{Frob_gen}. 
\item[{\rm (vi)}] If $d=2$, then $\F(S)=F_{2n}-F_{n+2}$ for the Fibonacci subsequence and $\F(S)=L_{2n+1}+L_{2n-1}-L_{n+2}$ for the Lucas subsequence; see Corollary \ref{Frob_special}. 
\item[{\rm (vii)}] If $d=2$, then $\g(S)$ for the Fibonacci subsequence is explicitly determined in terms of the solution of two recurrent sequences; see Theorem \ref{genus_special}. 
\end{itemize}
\vskip 10pt

\section{Preliminary Results} \label{prelim}
\vskip 10pt

A generalized Fibonacci sequence $\la V_n \ra_{n \ge 1}$ is defined by 
\begin{equation} \label{gen_def}
V_n = V_{n-1}+ V_{n-2}, \;\; n \ge 3, \quad \;\text{with}\; V_1 = a, V_2 = b, 
\end{equation}
where $a$ and $b$ are any positive integers. Two important special cases are (i) Fibonacci sequence $\{F_n\}_{n \ge 1}$ when $a=b=1$, and (ii) Lucas sequence $\{L_n\}_{n \ge 1}$ when $a=1$ and $b=3$. It is customary to extend these definitions to $F_0=F_2-F_1=0$ and $L_0=L_2-L_1=2$. The following formulae of Binet gives an explicit expression for both $F_n$ and $L_n$: 
\begin{equation} \label{binet}
F_n = \frac{{\phi}^n-{\ov{\phi}}^n}{\sqrt{5}}, \qquad L_n = {\phi}^n+{\ov{\phi}}^n, \quad n \ge 1, 
\end{equation}
where $\phi, \ov{\phi}$ are the zeros of $x^2-x-1$, $\phi>\ov{\phi}$. Since $|\ov{\phi}|<\frac{2}{3}$, $F_n$ is the integer closest to ${\phi}^n/\sqrt{5}$ and $L_n$ is the integer closest to ${\phi}^n$ when $n>1$. 
\vskip 5pt

\noindent Two of the most well known connections between the sequences $\{F_n\}$ and $\{L_n\}$ are: 
\[ L_n = F_{n+1}+F_{n-1}, \:n \ge 2, \quad \text{and} \quad F_{2n} = F_n L_n, \:n \ge 1. \]
The first may be easily derived from eqn.~\eqref{binet} or by induction, while the second is immediate from eqn.~\eqref{binet}. 
It is well known and easy to see, either by eqn.~\eqref{binet} or by induction, that 
\[ F_{n+1} F_{n-1} - F_n^2 = (-1)^n \quad \text{and} \quad L_{n+1} L_{n-1} - L_n^2 = 5(-1)^n, \]
for $n \ge 2$. 
\vskip 5pt

\noindent The following identities connecting generalized Fibonacci sequences with the Fibonacci sequence are useful in our subsequent work.  
\vskip 5pt

\begin{prop} \label{V_identity}
\begin{itemize}
\item[]
\item[{\rm (i)}]
For positive integers $m$ and $n$, 
\[ V_{m+n} = F_{n-1} V_m + F_n V_{m+1}. \]
In particular, $F_n \mid F_{kn}$ for each $k \ge 1$. 
\item[{\rm (ii)}]
For positive integers $k,n,d$, 
\[ F_d V_{n+kd} = (-1)^{d-1} F_{(k-1)d} V_n + F_{kd} V_{n+d}. \]
\item[{\rm (iii)}]
For positive integers $m$, $k$ and $d$, 
\[ V_{m+kd} \equiv V_{m+2}\,F_{kd} \!\!\!\!\pmod{V_m}. \]
\item[{\rm (iv)}]
For positive integers $n$ and $k$, 
\[ \sum_{i=1}^k V_{n+2i} = V_{n+2k+1} - V_{n+1}. \]
\end{itemize}
\end{prop}

\begin{Pf}
\begin{itemize}
\item[{\rm (i)}]
We fix $m$ and induct on $n$. The case $n=1$ is an identity and the case $n=2$ follows from the recurrence satisfied by the sequence $\{V_n\}$. Assuming the result for all positive integers less than $n$, we have 
\begin{eqnarray*} 
V_{m+n} & = & V_{m+(n-1)} + V_{m+(n-2)} \\
& = & \big(F_{n-2} V_m + F_{n-1} V_{m+1}\big) +  \big(F_{n-3} V_m + F_{n-2} V_{m+1}\big) \\
& = & \big(F_{n-2}+F_{n-3}\big)V_m + \big(F_{n-1}+F_{n-2}\big)V_{m+1} \\
& = & F_{n-1} V_m + F_n V_{m+1}. 
\end{eqnarray*}
This completes the proof by induction. 

In particular, with $V_n=F_n$ and $m=(k-1)n$, we have  
\[ F_{kn} = F_{(k-1)n}F_{n-1} + F_{(k-1)n+1}F_n. \]
So if $F_n \mid F_{(k-1)n}$, then $F_n \mid F_{kn}$. Hence, $F_n \mid F_{kn}$ for each $k \ge 1$ by induction. 

\item[{\rm (ii)}]
Interchanging $m$ and $n$ in part (i) and setting $m=d$ yields $F_d V_{n+1}=V_{n+d}-F_{d-1} V_n$. Therefore
\begin{eqnarray*}
F_d V_{n+kd} & = & F_d F_{kd-1} V_n + F_{kd} \big(V_{n+d}-F_{d-1} V_n\big) \\
& = & \big(F_d F_{kd-1} - F_{d-1} F_{kd}\big)V_n + F_{kd} V_{n+d} \\
& = & (-1)^{d-1} F_{(k-1)d} V_n + F_{kd} V_{n+d}. 
\end{eqnarray*}

\item[{\rm (iii)}]
This is obtained from part (i) by setting $n=kd$ since $V_{m+1}=V_{m+2}-V_m$. 

\item[{\rm (iv)}]
We fix $n$ and induct on $k$. The case $k=1$ follows from the recurrence satisfied by the sequence $\{V_n\}$. Assuming the result for all positive integers less than $k$, we have 
\begin{eqnarray*} 
\sum_{i=1}^k V_{n+2i} & = & \sum_{i=1}^{k-1} V_{n+2i} + V_{n+2k} \\
& = & \big( V_{n+2k-1} - V_{n+1} \big) + V_{n+2k} \\
& = & V_{n+2k+1} - V_{n+1}. 
\end{eqnarray*}
This completes the proof by induction. 
\end{itemize}
\end{Pf}
\vskip 5pt

\begin{cor}
If $\la V_n, V_{n+d}, V_{n+2d}, \ldots \ra$ is a numerical semigroup, then $\gcd(V_n,V_{n+d})=1$.  
\end{cor}

\begin{Pf}
The conclusion follows from the observation that $\gcd(V_n,V_{n+d})$ divides each term $V_{n+kd}$ by Proposition \ref{V_identity}, part (ii).    
\end{Pf}

\begin{cor} \label{V_gcd}
If $\gcd(V_1,V_2)=1$, then
\[ \gcd(V_n,V_{n+d}) = \gcd(V_n,F_d). \]
\end{cor}

\begin{Pf}
From the identity $V_{n+1}=V_n+V_{n-1}$ we have $\gcd(V_n,V_{n+1})=\gcd(V_{n-1},V_n)$. Since $\gcd(V_1,V_2)=1$, it follows that $\gcd(V_n,V_{n+1})=1$ for each positive integer $n$. Now Proposition \ref{V_identity}, part (i) gives $\gcd(V_n,V_{n+d})=\gcd(V_n,V_{n+1}F_d)=\gcd(V_n,F_d)$ since $\gcd(V_n,V_{n+1})=1$. 
\end{Pf}
\vskip 5pt

\noindent The following theorem provides a necessary and sufficient condition for $\la V_n, V_{n+d}, V_{n+2d}, \ldots \ra$ to form a numerical semigroup.  
\vskip 5pt

\begin{thm} \label{existence}
Let $S=\la V_n,V_{n+d},V_{n+2d}, \ldots \ra$. Then $S$ is a numerical semigroup if and only if $\gcd(V_1,V_2)=1$ and $\gcd(V_n,F_d)=1$.  
\end{thm}

\begin{Pf}
We recall that $S$ is a numerical semigroup if and only if $\gcd(V_n,V_{n+d},V_{n+2d},\ldots )=1$. Let $g=\gcd(a,b)$. Then $g$ divides each $V_m$ by the recurrence defining the terms in the sequence $\la V_m \ra$, and so a necessary condition for $S$ to be a numerical semigroup is $g=1$. 

Now assume $g=1$, and write $g^{\pr}=\gcd(V_n,F_d)$. By Proposition \ref{V_identity}, part (i), $V_{n+kd}=V_{n+(k-1)d}F_{d-1}+V_{n+(k-1)d+1}F_d$. Since $g^{\pr}$ divides $V_n$, if $g^{\pr}$ divides $V_{n+(k-1)d}$ then $g^{\pr}$ also divides $V_{n+kd}$. Hence, $g^{\pr}$ divides $V_{n+kd}$ for each $k \ge 0$ by induction, which implies $g^{\pr}=1$. This shows the necessity of the two $\gcd$ conditions. 

Conversely, $g^{\pr}=1$ and Corollary \ref{V_gcd} imply $\gcd(V_{n+d},V_n)=1$, so that $\gcd(V_n,V_{n+d},V_{n+2d},\ldots )=1$. This proves the sufficiency of the conditions.   
\end{Pf}
\vskip 10pt

\section{The Case where $d$ is odd} \label{odd}
\vskip 10pt

In this Section, we study the case where $d$ is odd. As a consequence of the identities in Proposition \ref{V_identity}, we show that each term $V_{n+kd}$, $k \ge 2$, is of the form $V_n x+V_{n+d} y$ with $x, y \in {\Z}_{\ge 0}$. This enables us to easily determine $\F(S)$, $\g(S)$ and the set $\PF(S)$ in this case. 
\vskip 5pt
 
\begin{thm} \label{e(S)_odd}
If $\gcd(V_1,V_2)=1$ and $\gcd(V_n,F_d)=1$, then 
\[  \la V_n,V_{n+d},V_{n+2d}, \ldots \ra = \la V_n, V_{n+d} \ra. \]
\end{thm}
   
\begin{Pf}
The two gcd conditions are necessary and sufficient to ensure $S=\la V_n,V_{n+d},V_{n+2d}, \ldots \ra$ is a numerical semigroup by Theorem \ref{existence}. In order to prove the result of this theorem, we must only show that $V_{n+kd} \in \la V_n, V_{n+d} \ra$ for each $k \ge 2$. This is a direct consequence of Proposition \ref{V_identity}, part (ii), since $F_d$ divides both $F_{(k-1)d}$ and $F_{kd}$ for each $k \ge 2$. 
\end{Pf}

\begin{thm} \label{odd_case}
Let $d$ be odd, $\gcd(V_1,V_2)=1$ and $\gcd(V_n,F_d)=1$. If $S=\la V_n,V_{n+d},V_{n+2d}, \ldots \ra$, then
\begin{itemize}
\item[{\rm (i)}]
\[ \e(S) = 2. \]
\item[{\rm (ii)}]
\[ \F(S) = (V_n-1)(V_{n+d}-1)-1. \]
\item[{\rm (iii)}]
\[ \g(S) = \frac{1}{2}(V_n-1)(V_{n+d}-1).  \]
\item[{\rm (iv)}]
\[ \PF(S) = \Big\{ (V_n-1)(V_{n+d}-1)-1 \Big\}. \]
\end{itemize}
\end{thm}

\begin{Pf}
Part (i) is the result in Theorem \ref{e(S)_odd}. Parts (ii), (iii) and (iv) are consequences of part (i) and eqn.~\eqref{e(S)=2}. 
\end{Pf}
\vskip 10pt

\section{The Case where $d=2$} \label{d=2}
\vskip 10pt

In this Section, we study the case where $d=2$. This case turns out to be a lot more challenging than when $d$ is odd, in part due to the fact that the crucial identity that allows each $V_{n+kd} \in \la V_n, V_{n+d} \ra$ depends on the parity of $d$. The following result leads to the determination of the Ap\'{e}ry set for $S$ with respect to $V_n$. 
\vskip 5pt
 
\begin{prop} \label{2-representation}
Fix $x \in \{1,\ldots,V_n-1\}$, and let $k$ be such that $F_{2k} \le x < F_{2k+2}$.  Then there exists ${\l}_1,\ldots,{\l}_k$, with each ${\l}_i \in \{0,1,2\}$ and ${\l}_k \ge 1$ such that 
\[ x = \sum_{i=1}^k {\l}_i F_{2i}, \quad s = \sum_{i=1}^k {\l}_i V_{n+2i}, \quad V_{n+2k} \le s < V_{n+2k+2}, \quad s \equiv V_{n+2}\,x\!\!\!\!\pmod{V_n}. \] 
\end{prop}

\begin{Pf}
We define the sequence ${\l}_k, {\l}_{k-1}, \ldots, {\l}_1$ by using the Greedy Algorithm on $x$ with respect to the sequence $F_2, F_4, F_6, \ldots, F_{2k}$: 
\begin{equation} \label{greedy}
{\l}_k = \left\lf \frac{x}{F_{2k}} \right\rf, \quad {\l}_j = \left\lf \frac{x-\sum_{i=j+1}^k {\l}_i F_{2i}}{F_{2j}} \right\rf, \:\: j=k-1,k-2,\ldots,1. 
\end{equation}
Since $F_{2k+2}=(F_{2k}+F_{2k-1})+F_{2k}$ and $F_{2k} \le x < F_{2k+2}$, we have $1 \le {\l}_k \le 2$. For each $j \in \{1,\ldots,k-1\}$, since 
\[ {\l}_{j+1} F_{2j+2} \le x - \sum_{i=j+2}^k {\l}_i F_{2i} < ({\l}_{j+1}+1) F_{2j+2}, \]
we have
\[ 0 \le x - \sum_{i=j+1}^k {\l}_i F_{2i} = x - \sum_{i=j+2}^k {\l}_i F_{2i} - {\l}_{j+1} F_{2j+2} < F_{2j+2}, \]
so that 
\[ 0 \le \frac{x-\sum_{i=j+1}^k {\l}_i F_{2i}}{F_{2j}} < \frac{F_{2j+2}}{F_{2j}} = 2 + \frac{F_{2j-1}}{F_{2j}}. \]
Thus, ${\l}_j \in \{0,1,2\}$. 

Now 
\[ {\l}_1 = \left\lf \frac{x-\sum_{i=2}^k {\l}_i F_{2i}}{F_2} \right\rf = x-\sum_{i=2}^k {\l}_i F_{2i}, \]
so $x=\sum_{i=1}^k {\l}_i F_{2i}$. 

Define $s=\sum_{i=1}^k {\l}_i V_{n+2i}$. By Proposition \ref{V_identity}, part (ii), 
\[ s = \sum_{i=1}^k {\l}_i V_{n+2i} \equiv \sum_{i=1}^k {\l}_i V_{n+2}\,F_{2i} = V_{n+2}\,x \!\!\!\!\pmod{V_n}. \] 
Since ${\l}_k \ge 1$, we have $s \ge V_{n+2k}$. To prove the upper bound for $s$, we consider two cases: (i) ${\l}_k=1$, and (ii) ${\l}_k=2$. 
\vskip 5pt

\noindent {\sc Case} (I): If ${\l}_k=1$, then 
\[ s \le 2 \sum_{i=1}^{k-1} V_{n+2i} + V_{n+2k} \le 2 \big( V_{n+2k-1} - V_{n+1} \big) + V_{n+2k} < \big( V_{n+2k}+V_{n+2k-1} \big) + V_{n+2k-1} < V_{n+2k+2} \] 
using Proposition \ref{V_identity}, part (iii). 
\vskip 5pt

\noindent {\sc Case} (II): Suppose ${\l}_k=2$. We claim that one of the following cases must arise: (i) ${\l}_i=1$ for $i \in \{1,\ldots,k-1\}$; (ii) there exists $r \in \{1,\ldots,k-1\}$ such that ${\l}_r=0$ and ${\l}_i=1$ for $i \in \{r+1,\ldots,k-1\}$. 

\noindent If neither of these cases is true, then there must exist $t \in \{1,\ldots,k-1\}$ such that ${\l}_t=2$ and ${\l}_i=1$ for $i \in \{t+1,\ldots,k-1\}$. But then 
\[ x = \sum_{i=1}^k {\l}_i F_{2i} \ge  \sum_{i=t}^k F_{2i} + F_{2t} + F_{2k} = \big( F_{2k+1} - F_{2t-1} \big) + F_{2t} + F_{2k} \ge F_{2k+2} \]
using Proposition \ref{V_identity}, part (iii) with $V_n=F_n$. This contradiction proves the claim. 
\vskip 5pt

\noindent In case (i), we have 
\[ s = \sum_{i=1}^k {\l}_i V_{n+2i} = V_{n+2k} + \sum_{i=1}^k V_{n+2i} = V_{n+2k} + \big( V_{n+2k-1} - V_{n+1} \big) < V_{n+2k+2} \]
using Proposition \ref{V_identity}, part (iii). 

\noindent In case (ii), 
\begin{eqnarray*} 
s & = & \sum_{i=1}^k {\l}_i V_{n+2i} \\
& \le & V_{n+2k} + \sum_{i=r+1}^k V_{n+2i} + 2 \sum_{i=1}^{r-1} V_{n+2i} \\
& = & V_{n+2k} + \big( V_{n+2k+1} - V_{n+2r+1} \big) + 2 \big( V_{n+2r-1} - V_{n+1} \big) \\
& = & V_{n+2k+2} - \left( \big(V_{n+2r+1} - V_{n+2r-1}\big) - V_{n+2r-1} \right) - 2V_{n+1} \\
& < & V_{n+2k+2}.  
\end{eqnarray*}

\noindent This completes the proof of the Proposition. 
\end{Pf}
\vskip 5pt

\begin{notn} \label{s(x)} 
The sequence ${\l}_1, \ldots, {\l}_k$ in the proof of Proposition \ref{2-representation} is determined by applying the Greedy Algorithm to $x$ with respect to the sequence $F_2, F_4, F_6, \ldots, F_{2k}$, and $s$ is then determined from this sequence. So if $x=\sum_{i=1}^k {\l}_i F_{2i}$, then $s=\sum_{i=1}^k {\l}_i V_{n+2i}$. We use the expression $s(x)$ to show the dependence of $s$ on $x$ via the sequence ${\l}_1,\ldots,{\l}_k$.   
\end{notn}
\vskip 5pt

\begin{lem} \label{s_special}
We have $s(1)=V_{n+2}$ and $s(2)=2V_{n+2}$. For each positive integer $m$, we have
\begin{itemize}
\item[{\rm (i)}]
\[ s(F_m) = \begin{cases} 
                    V_{n+m} + V_n  & \mbox{ if $m$ is odd}; \\ 
                    V_{n+m} & \mbox{ if $m$ is even}. 
                  \end{cases}  
\]
\item[{\rm (ii)}]
\[ s(L_m) = \begin{cases} 
                    V_{n+m+1} + V_{n+m-1} & \mbox{ if $m$ is odd, $m>1$}; \\ 
                    V_{n+m+1} + V_{n+m-1} + V_n & \mbox{ if $m$ is even}.  
                   \end{cases}  
\]
\item[{\rm (iii)}]
\[ s(F_m-1) = V_{n+m} - V_{n+1}, m>2. \] 
\item[{\rm (iv)}]
\[ s(L_m-1) = V_{n+m+1} + V_{n+m-1} - V_{n+1}, m>2. \] 
\end{itemize}
\end{lem}

\begin{Pf}
It is easy to verify that $s(1)=V_{n+2}$ and $s(2)=2V_{n+2}$. 
\begin{itemize}
\item[{\rm (i)}]
If $m=2k+1$, then $F_{2k} \le F_{2k+1}<F_{2k+2}$, and from eqn.~\eqref{greedy}, ${\l}_k=\left\lf \frac{F_{2k+1}}{F_{2k}} \right\rf=1$. Now ${\l}_{k-1}=\left\lf \frac{F_{2k+1}-F_{2k}}{F_{2k-2}} \right\rf=\left\lf\frac{F_{2k-1}}{F_{2k-2}}\right\rf=1$. 

If ${\l}_i=\left\lf\frac{F_{2i+1}}{F_{2i}}\right\rf$ for $i=k,k-1,\ldots,j+1$, then ${\l}_j=\left\lf\frac{F_{2j+3}-F_{2j+2}}{F_{2j}}\right\rf=\left\lf\frac{F_{2j+1}}{F_{2j}}\right\rf=1$. Finally, ${\l}_1=\left\lf\frac{F_3}{F_2}\right\rf=2$. 

Thus, $s(F_{2k+1})=\sum_{i=1}^k V_{n+2i}+V_{n+2}=\big(V_{n+2k+1}-V_{n+1}\big)+V_{n+2}=V_{n+2k+1}+V_n$ by Proposition \ref{V_identity}, part (iii).   

If $m=2k$, then $F_{2k}=F_m<F_{2k+2}$, and from eqn.~\eqref{greedy}, ${\l}_k=1$ and ${\l}_i=0$ for $i<k$. Thus, $s(F_{2k})=V_{n+2k}$.   

\item[{\rm (ii)}]
If $m=2k-1$, then $L_m=F_{2k}+F_{2k-2}$, so that $s(L_{2k-1})=V_{n+2k}+V_{n+2k-2}$.   

If $m=2k$, then $F_{2k}<L_{2k}=F_{2k+1}+F_{2k-1}<F_{2k+2}$. Since $F_{2k+1}+F_{2k-1}=F_{2k}+\big(F_{2k-1}+F_{2k-2}\big)+\big(F_{2k-1}-F_{2k-2}\big)=2F_{2k}+F_{2k-3}$, from eqn.~\eqref{greedy}, ${\l}_k=\left\lf \frac{L_{2k}}{F_{2k}} \right\rf=2$, and ${\l}_{k-1}=\left\lf \frac{L_{2k}-2F_{2k}}{F_{2k-2}} \right\rf=\left\lf \frac{F_{2k-3}}{F_{2k-2}} \right\rf=0$. Now ${\l}_i=1$ for $i=k-2,\ldots,2$ and ${\l}_1=2$ by the odd case in part (i). 

Thus, $s(L_{2k})=\sum_{i=1}^k V_{n+2i}+V_{n+2k}-V_{n+2k-2}+V_{n+2}=\big(V_{n+2k+1}-V_{n+1}\big)+V_{n+2k-1}+V_{n+2}=V_{n+2k+1}+V_{n+2k-1}+V_n$ by Proposition \ref{V_identity}, part (iii).   

\item[{\rm (iii)}]
If $m=2k+1$, then $F_{2k} \le  F_m-1<F_{2k+1}<2F_{2k}<F_{2k+2}$, and from eqn.~\eqref{greedy}, ${\l}_k=\left\lf \frac{F_{2k+1}-1}{F_{2k}} \right\rf=1$. Now ${\l}_{k-1}=\left\lf \frac{(F_{2k+1}-1)-F_{2k}}{F_{2k-2}} \right\rf=\left\lf\frac{F_{2k-1}-1}{F_{2k-2}}\right\rf=1$, by the argument for ${\l}_k$. 

If ${\l}_i=\left\lf\frac{F_{2i+1}-1}{F_{2i}}\right\rf$ for $i=k,k-1,\ldots,j+1$, then ${\l}_j=\left\lf\frac{(F_{2j+3}-1)-F_{2j+2}}{F_{2j}}\right\rf=\left\lf\frac{F_{2j+1}-1}{F_{2j}}\right\rf=1$. Finally, ${\l}_1=\left\lf\frac{F_3-1}{F_2}\right\rf=1$. 

Thus, $s(F_{2k+1}-1)=\sum_{i=1}^k V_{n+2i}=V_{n+2k+1}-V_{n+1}$ by Proposition \ref{V_identity}, part (iii).   

If $m=2k+2$, then $F_{2k}<2F_{2k} \le F_m-1=F_{2k+2}-1<F_{2k+2}$, and from eqn.~\eqref{greedy}, ${\l}_k=\left\lf \frac{F_{2k+2}-1}{F_{2k}} \right\rf=2$. Now ${\l}_{k-1}=\left\lf \frac{(F_{2k+2}-1)-2F_{2k}}{F_{2k-2}} \right\rf=\left\lf\frac{F_{2k-1}-1}{F_{2k-2}}\right\rf=1$, by the argument for the odd case. Arguing as in the odd case, each ${\l}_i=1$ for $i=k-1,k-2,\ldots,1$. 

Thus, $s(F_{2k+2}-1)=\sum_{i=1}^k V_{n+2i}+V_{n+2k}=\left(V_{n+2k+1}-V_{n+1}\right)+V_{n+2k}=V_{n+2k+2}-V_{n+1}$ by Proposition \ref{V_identity}, part (iii).   

\item[{\rm (iv)}]
If $m=2k-1$, then $F_{2k} \le L_m-1=F_{2k}+F_{2k-2}-1 < F_{2k+2}$, and from eqn.~\eqref{greedy}, ${\l}_k=\left\lf \frac{L_{2k-1}-1}{F_{2k}} \right\rf=1$. Now ${\l}_{k-1}=\left\lf \frac{L_{2k-1}-1-F_{2k}}{F_{2k-2}} \right\rf=\left\lf\frac{F_{2k-2}-1}{F_{2k-2}}\right\rf=0$, so that ${\l}_{k-2}=\left\lf \frac{F_{2k-2}-1}{F_{2k-4}} \right\rf=\left\lf\frac{2F_{2k-4}+F_{2k-5}-1}{F_{2k-4}}\right\rf=2$. Hence ${\l}_{k-3}=\left\lf\frac{F_{2k-5}-1}{F_{2k-6}}\right\rf$, and so ${\l}_i=1$ for $i=k-3,\ldots,1$ by the odd case in part (iii). 

Thus, $s(L_{2k-1}-1)=\sum_{i=1}^k V_{n+2i}-V_{n+2k-2}+V_{n+2k-4}=\left(V_{n+2k+1}-V_{n+1}\right)-V_{n+2k-3}$ by Proposition \ref{V_identity}, part (iii). So, for odd $m$, $s(L_m-1)=V_{n+m+2}-V_{n+m-2}-V_{n+1}=V_{n+m+1}+\big(V_{n+m}-V_{n+m-2}\big)-V_{n+1}=V_{n+m+1}+V_{n+m-1}-V_{n+1}$.   

If $m=2k$, then $F_{2k}<L_{2k}-1=F_{2k+1}+F_{2k-1}-1<F_{2k+2}$. Since $F_{2k+1}+F_{2k-1}-1=F_{2k}+\big(F_{2k-1}+F_{2k-2}\big)+\big(F_{2k-1}-F_{2k-2}-1\big)=2F_{2k}+F_{2k-3}-1$, from eqn.~\eqref{greedy}, ${\l}_k=\left\lf \frac{L_{2k}-1}{F_{2k}} \right\rf=2$, and ${\l}_{k-1}=\left\lf \frac{L_{2k}-1-2F_{2k}}{F_{2k-2}} \right\rf=\left\lf \frac{F_{2k-3}-1}{F_{2k-2}} \right\rf=0$. Now ${\l}_i=1$ for $i=k-2,\ldots,2$ and ${\l}_1=2$ by the odd case in part (i). 

Thus, $s(L_{2k}-1)=\sum_{i=1}^k V_{n+2i}+V_{n+2k}-V_{n+2k-2}=\big(V_{n+2k+1}-V_{n+1}\big)+V_{n+2k-1}$ by Proposition \ref{V_identity}, part (iii), and again $s(L_m-1)=V_{n+m+1}+V_{n+m-1}-V_{n+1}$ when $m$ is even.    
\end{itemize}
\end{Pf}
\vskip 5pt

\noindent The following result proves that the Greedy Algorithm for an arbitrary positive integer $x$ with respect to the sequence $F_2, F_4, F_6, \ldots $ employed to compute $s(x)$ is optimal. 
\vskip 5pt

\begin{thm} \label{Greedy_is_optimum}
For any sequence ${\a}_1, \ldots, {\a}_m$ of nonnegative integers, not all zero, 
\[ s\left( \sum_{i=1}^m {\a}_i F_{2i} \right) \le \sum_{i=1}^m {\a}_i V_{n+2i}. \]
\end{thm}

\begin{Pf}
We induct on the sum $\s=\sum_{i=1}^m {\a}_i F_{2i}$. If $\s=1$, then $m={\a}_1=1$ and the two sides are equal. For some positive integer $\s$, assume the result holds whenever the sum $\sum_{i=1}^m {\a}_i F_{2i}<\s$. 
\vskip 5pt

\noindent Using Proposition \ref{2-representation}, define the sequence ${\l}_1,\ldots,{\l}_k$ for $x=\s$. Suppose ${\a}_1,\ldots,{\a}_m$ is any sequence of nonegative integers such that $\s=\sum_{i=1}^m {\a}_i F_{2i}$; we may assume that ${\a}_m \ge 1$. Note that $m \le k$, for if $m>k$, then $\sum_{i=1}^m {\a}_i F_{2i} \ge F_{2k+2}>\s$. 
\vskip 5pt

\noindent If $m=k$, then $1 \le {\a}_k \le {\l}_k$. By Induction Hypothesis, 
\[ s\left( \sum_{i=1}^k {\a}_i F_{2i} - F_{2k} \right) \le \sum_{i=1}^k {\a}_i V_{n+2i} - V_{n+2k}. \]
Since  
\[ s\left( \sum_{i=1}^k {\a}_i F_{2i} - F_{2k} \right) =  s\left( \sum_{i=1}^k {\l}_i F_{2i} - F_{2k} \right) = \sum_{i=1}^k {\l}_i V_{n+2i} - V_{n+2k} = s\left( \sum_{i=1}^k {\l}_i F_{2i} \right) - V_{n+2k}, \]
we have 
\[ s\left( \sum_{i=1}^k {\a}_i F_{2i} \right) = s\left( \sum_{i=1}^k {\l}_i F_{2i} \right) \le \sum_{i=1}^k {\a}_i V_{n+2i}. \]
This proves the Proposition when $m=k$. 
\vskip 5pt

\noindent Suppose $m<k$. By Induction Hypothesis, 
\begin{equation} \label{IH}
s\left( \sum_{i=1}^m {\a}_i F_{2i} - F_{2m} \right) \le \sum_{i=1}^m {\a}_i V_{n+2i} - V_{n+2m}. 
\end{equation}
Two cases arise: (I) $\sum_{i=1}^m {\a}_i F_{2i} - F_{2m} \ge F_{2k}$, and (II) $\sum_{i=1}^m {\a}_i F_{2i} - F_{2m} < F_{2k}$. 
\vskip 5pt

\noindent {\sc Case} (I): Suppose $\s-F_{2m}=\sum_{i=1}^m {\a}_i F_{2i} - F_{2m} \ge F_{2k}$. Let ${\l}_1^{\pr}, \ldots, {\l}_k^{\pr}$ be the sequence determined by the Greedy Algorithm for $\s-F_{2m}$. Then 
\begin{equation} \label{m<k_subcase1}
s\left( \sum_{i=1}^m {\a}_i F_{2i} - F_{2m} \right) = s\left( \sum_{i=1}^k {\l}_i^{\pr} F_{2i} \right) = \sum_{i=1}^k {\l}_i^{\pr} V_{n+2i}. 
\end{equation}
If we replace ${\l}_m^{\pr}$ by ${\l}_m^{\pr}+1$ and retain the other ${\l}_i^{\pr}$, and apply the case $m=k$ discussed above, we get 
\[ s\left( \sum_{i=1}^m {\a}_i F_{2i} \right) =  s\left( \sum_{i=1}^k {\l}_i^{\pr} F_{2i} + F_{2m} \right) \le \sum_{i=1}^k {\l}_i^{\pr} V_{n+2i} + V_{n+2m} \le \sum_{i=1}^m {\a}_i V_{n+2i}. \]
from eqn.~\eqref{IH} and eqn.~\eqref{m<k_subcase1}. This proves Case (i). 
\vskip 5pt

\noindent {\sc Case} (II): Suppose $\s-F_{2m}=\sum_{i=1}^m {\a}_i F_{2i} - F_{2m} < F_{2k}$. Since $\s-F_{2m} \ge F_{2k}-F_{2k-2}>F_{2k-2}$, the sequence determined by the Greedy Algorithm for $\s-F_{2m}$ is ${\l}_1^{\pr}, \ldots, {\l}_{k-1}^{\pr}$. 
\vskip 5pt

\noindent Note that $\s-F_{2k-2}$ lies between $F_{2k}-F_{2k-2}=F_{2k-1}$ and $F_{2k}$. Let ${\l}_1^{\pr\pr}, \ldots, {\l}_{k-1}^{\pr\pr}$ be the sequence determined by the Greedy Algorithm for $\s-F_{2k-2}$. We claim that one of the following cases must arise: (i) ${\l}_i^{\pr\pr}=1$ for $i \in \{1,\ldots,k-1\}$; (ii) there exists $r \in \{1,\ldots,k-1\}$ such that ${\l}_r^{\pr\pr}=2$ and ${\l}_i^{\pr\pr}=1$ for $i \in \{r+1,\ldots,k-1\}$. 
\vskip 5pt

\noindent If neither of these cases is true, then there must exist $t \in \{1,\ldots,k-1\}$ such that ${\l}_t^{\pr\pr}=0$ and ${\l}_i^{\pr\pr}=1$ for $i \in \{t+1,\ldots,k-1\}$. But then 
\begin{eqnarray*} 
{\l}_t^{\pr\pr} & = & \left\lf \frac{\s-F_{2k-2}-\sum_{i=t+1}^{k-1} {\l}_i^{\pr\pr} F_{2i}}{F_{2t}} \right\rf \ge \left\lf \frac{F_{2k}-F_{2k-2}-\sum_{i=t+1}^{k-1} F_{2i}}{F_{2t}} \right\rf \\
& \ge & \left\lf \frac{F_{2k}-F_{2k-2}-\big(F_{2k-1}-F_{2t+1}\big)}{F_{2t}} \right\rf = 1. 
\end{eqnarray*}
This contradiction proves the claim. 
\vskip 5pt

\noindent In Case (i), using Proposition \ref{V_identity}, part (iii), we have 
\[ \sum_{i=1}^{k-1} {\l}_i^{\pr\pr} F_{2i} = F_{2k-1} - F_1 < F_{2k} - F_{2k-2}, \]
contradicting the fact that ${\l}_1^{\pr\pr}, \ldots, {\l}_{k-1}^{\pr\pr}$ is the sequence determined by the Greedy Algorithm for $\s-F_{2k-2}$ and $\s \ge F_{2k}$. This rules out Case (i).  
\vskip 5pt

\noindent In Case (ii), using Proposition \ref{V_identity}, part (iii), we get
\begin{eqnarray} \label{m<k_subcase2}
s\left( \sum_{i=1}^{k-1} {\l}_i^{\pr\pr} F_{2i} \right) + V_{n+2k-2} & = & \sum_{i=1}^{k-1} {\l}_i^{\pr\pr} V_{n+2i} + V_{n+2k-2} \nonumber \\
& = & \sum_{i=1}^{r-1} {\l}_i^{\pr\pr} V_{n+2i} + \sum_{i=r+1}^{k-2} V_{n+2i} + V_{n+2k-2} + 2V_{n+2r} + V_{n+2k-2} \nonumber \\
& = & \sum_{i=1}^{r-1} {\l}_i^{\pr\pr} V_{n+2i} + \big( V_{n+2k-3} - V_{n+2r+1} \big) + V_{n+2k-2} + 2V_{n+2r} + V_{n+2k-2} \nonumber \\
& = & \sum_{i=1}^{r-1} {\l}_i^{\pr\pr} V_{n+2i} + V_{n+2k} + V_{n+2r-2}. 
\end{eqnarray}
\vskip 5pt 

\noindent We have 
\begin{eqnarray*}
\s - F_{2k-2} & = & \sum_{i=1}^{k-1} {\l}_i^{\pr\pr} F_{2i} \\
& = & \sum_{i=1}^{r-1} {\l}_i^{\pr\pr} F_{2i} + \sum_{i=r+1}^{k-1} F_{2i} + 2F_{2r} \\
& = & \sum_{i=1}^{r-1} {\l}_i^{\pr\pr} F_{2i} + \big( F_{2k-1} - F_{2r+1} \big) + 2F_{2r} \\
& = & \sum_{i=1}^{r-1} {\l}_i^{\pr\pr} F_{2i} + \big( F_{2k} - F_{2k-2} \big) + F_{2r-2}. 
\end{eqnarray*}
\vskip 5pt

\noindent By the Induction Hypothesis, 
\[ s\left( \sum_{i=1}^{k-1} {\l}_i^{\pr} F_{2i} - F_{2k-2} + F_{2m} \right) \le \sum_{i=1}^{k-1} {\l}_i^{\pr} V_{n+2i} - V_{n+2k-2} + V_{n+2m}. \]
Applying the case $m=k$ discussed above to $s\left(\s\right)$ and using eqn.~\eqref{m<k_subcase2}, we have  
\begin{eqnarray*} 
s\left( \sum_{i=1}^{r-1} {\l}_i^{\pr\pr} F_{2i} + F_{2k} + F_{2r-2} \right) & \le & \sum_{i=1}^{r-1} {\l}_i^{\pr\pr} V_{n+2i} + V_{n+2k} + V_{n+2r-2} \\
& = & s\left( \sum_{i=1}^{k-1} {\l}_i^{\pr\pr} F_{2i} \right) + V_{n+2k-2} \\
& = & s\left( \sum_{i=1}^{k-1} {\l}_i^{\pr} F_{2i} - F_{2k-2} + F_{2m} \right) + V_{n+2k-2} \\
& \le & \sum_{i=1}^{k-1} {\l}_i^{\pr} V_{n+2i} + V_{n+2m} \\
& = & s\left( \sum_{i=1}^m {\a}_i F_{2i} - F_{2m} \right) + V_{n+2m} \\
& \le & \sum_{i=1}^m {\a}_i V_{n+2i}.   
\end{eqnarray*} 

\noindent This completes Case (ii), and the proof. 
\end{Pf}
\vskip 5pt

\begin{lem} \label{s_increasing}
For any positive integer $m$, $s(m)<s(m+1)$. 
\end{lem}

\begin{Pf}
We induct on $m$. By Proposition \ref{s_special}, part (i), $V_{n+2}=s(1)<s(2)=2V_{n+2}$. Assume $s(i-1)<s(i)$ for $1 \le i \le m$. If $m=F_{2k}-1$ for some $k$, then $s(m)=V_{n+2k}-V_{n+1} <V_{n+2k}=s(m+1)$ by Proposition \ref{s_special}, part (i). Otherwise $F_{2k} \le m<F_{2k+2}$, and so $s(m)=s(m-F_{2k})+V_{n+2k}$ while $s(m+1)=s(m+1-F_{2k})+V_{n+2k}$. By Induction Hypothesis, $s(m-F_{2k})<s(m+1-F_{2k})$, so that $s(m)<s(m+1)$, proving the Proposition by induction. 
\end{Pf}
\vskip 5pt

\begin{thm} \label{Apery_set}
Let $\gcd(V_1,V_2)=1$. The Ap\'{e}ry set for $S=\la V_n, V_{n+2}, V_{n+4}, \ldots \ra$ is given by  
\[ \Ap(S,V_n) = \{ s(x): 1 \le x \le V_n-1 \} \cup \{0\}. \] 
\end{thm}

\begin{Pf}
For $x \in \{1,\ldots,V_n-1\}$, we show that $s(x)$ is the least positive integer in $S$ that is congruent to $V_{n+2}\,x$ modulo $V_n$. This proves the result since $\{V_{n+2}\,x: 1 \le x \le V_n-1\}$ is the set of non-zero residues modulo $V_n$ as $\gcd(V_n,V_{n+2})=1$. 

Suppose $s \in S$ is congruent to $V_{n+2}\,x$ modulo $V_n$. Then $s=\sum_{i \ge 0} {\a}_i V_{n+2i}$, with each ${\a}_i \ge 0$. Since $s \equiv V_{n+2}\,x\pmod{V_n}$, we have $\sum_{i \ge 1} {\a}_i F_{2i} \equiv x\pmod{V_n}$ as $\gcd(V_n,V_{n+2})=1$. Since $x \le V_n-1$, we have $x \le \sum_{i \ge 1} {\a}_i F_{2i}$, so that  
\[ s(x) \le s\left( \sum_{i \ge 1} {\a}_i F_{2i}\right) \le \sum_{i \ge 1} {\a}_i V_{n+2i} \le s \]
by Theorem \ref{Greedy_is_optimum} and Lemma \ref{s_increasing}. 
\end{Pf}
\vskip 5pt

\begin{thm} \label{Frob_gen}
Let $\gcd(V_1,V_2)=1$. If $S=\la V_n, V_{n+2}, V_{n+4}, \ldots \ra$, then 
\begin{itemize}
\item[{\rm (i)}]
\[ \F(S) = s(V_n-1) - V_n, \]
\item[{\rm (ii)}]
\[ \g(S) = \frac{1}{V_n} \left( \sum_{x=1}^{V_n-1} s(x) \right) - \frac{V_n-1}{2}. \] 
\end{itemize}
\end{thm}

\begin{Pf}
These are direct consequences of Proposition \ref{prelim}, Theorem \ref{Apery_set} and Lemma \ref{s_increasing}. 
\end{Pf}
\vskip 5pt
 
\newpage 

\begin{cor} \label{Frob_special}
\begin{itemize}
\item[]
\item[{\rm (i)}]
If $S=\la F_n, F_{n+2}, F_{n+4}, \ldots \ra$, $n \ge 3$, then $\F(S)=F_{2n}-F_{n+2}$.
\item[{\rm (ii)}]
If $S=\la L_n, L_{n+2}, L_{n+4}, \ldots \ra$, $n \ge 4$, then $\F(S)=L_{2n+1}+L_{2n-1}-L_{n+2}$. 
\end{itemize}
\end{cor}

\begin{Pf}
This is a direct consequence of Theorem \ref{Frob_gen} and Proposition \ref{s_special}.  
\end{Pf}
\vskip 5pt

\noindent We determine the embedding dimension for $S$ by using the characterization given in the Introduction. We use the Greedy Algorithm to show that every element in $\Ap(S,V_n)$ can be expressed as a nonnegative linear combination of elements of the set which we claim is the minimal generating set for $S$. 
\vskip 5pt

\begin{thm} \label{embed_dim}
If $\gcd(V_1,V_2)=1$, then the minimal set of generators for $S=\la V_n, V_{n+2}, V_{n+4},\ldots \ra$ is 
\[  A = \{ V_n, V_{n+2}, \ldots, V_{n+2(\k-1)} \}, \]
where $\k$ is given by $F_{2(\k-1)} \le V_n-1 < F_{2\k}$.  
\end{thm}

\begin{Pf}
Let $\k$ be given by $F_{2(\k-1)} \le V_n-1 < F_{2\k}$. We show that $A=\{ V_n,V_{n+2}, \ldots,V_{n+2(\k-1)}\}$ is the minimal set of generators for $S$ by using the result in Proposition \ref{prelim}, part (iii). 
\vskip 5pt

\noindent We first show that every element in $\Ap(S,V_n)$ can be expressed as a nonnegative linear combination of elements of $A$. 
\vskip 5pt

\noindent Fix $x \in \{1,\ldots,V_n-1\}$, so that $x<F_{2\k}$. Applying the Greedy Algorithm on $x$ with respect to $F_2, F_4, F_6, \ldots $ gives $x = \sum_{i=1}^k {\l}_i F_{2i}$ and $s(x) = \sum_{i=1}^k {\l}_i V_{n+2i}$, with $k<\k$. Since $\Ap(S,V_n)=\{s(x): 1 \le x \le V_n-1\}\:\cup\:\{0\}$, we have the claim. 
\vskip 5pt

\noindent To show that $A$ is a minimal set of generators for $S$, we show that $V_{n+2k} \notin \la V_n, V_{n+2}, \ldots, V_{n+2(k-1)} \ra$ for each $k \in \{1,\ldots,\k-1\}$. Suppose there exist nonnegative integers ${\mu}_0,\ldots,{\mu}_{k-1}$ such that 
\begin{equation} \label{embed_dim1}
{\mu}_0 V_n + {\mu}_1 V_{n+2} + \cdots + {\mu}_{k-1} V_{n+2(k-1)} = \sum_{i=0}^{k-1} {\mu}_i V_{n+2i} = V_{n+2k}. 
\end{equation}
Reducing both sides of eqn.~\eqref{embed_dim1} modulo $V_n$ and using Proposition \ref{V_identity}, part (ii) gives
\[  V_{n+2} \sum_{i=1}^{k-1} {\mu}_i F_{2i} \equiv F_{2k} V_{n+2} \!\!\!\!\pmod{V_n}. \]
Since $\gcd(V_n,V_{n+2})=1$, we have 
\[ \sum_{i=1}^{k-1} {\mu}_i F_{2i} \equiv F_{2k} \!\!\!\!\pmod{V_n}. \]
Since $k \le \k-1$, we have $F_{2k} \le F_{2(\k-1)} < V_n$, and so $\sum_{i=1}^{k-1} {\mu}_i F_{2i} \ge F_{2k}$. By Theorem \ref{Greedy_is_optimum}, Proposition \ref{s_increasing} and eqn.~\eqref{embed_dim1}, we have 
\[ V_{n+2k} = s\left( F_{2k} \right) \le s\left( \sum_{i=1}^{k-1} {\mu}_i F_{2i} \right) \le  \sum_{i=0}^{k-1} {\mu}_i V_{n+2i} = V_{n+2k}. \]  
Therefore, by Proposition \ref{s_increasing}, 
\[ F_{2k} = \sum_{i=1}^{k-1} {\mu}_i F_{2i} = \sum_{i=1}^m {\mu}_i F_{2i}, \]
where $m$ be the largest integer for which ${\mu}_m>0$. Then by Theorem \ref{Greedy_is_optimum}
\[ s\left( F_{2k}-F_{2m} \right) \le \sum_{i=1}^m {\mu}_i V_{n+2i} - V_{n+2m}. \]
\vskip 5pt

\noindent Since $F_{2k}-F_{2m} \ge F_{2k}-F_{2k-2} > F_{2k-2}$, the Greedy Algorithm applied to $F_{2k}-F_{2m}$ yields a sequence ${\l}_1,\ldots,{\l}_{k-1}$, with ${\l}_{k-1}>0$. Applying Theorem \ref{Greedy_is_optimum} to the sequence obtained by adding $1$ to ${\l}_m$ and subtracting $1$ from ${\l}_{k-1}$ yields  
\[  s\left( F_{2k}-F_{2k-2} \right) \le s\left( F_{2k}-F_{2m} \right) + V_{n+2m} - V_{n+2k-2}. \] 
\vskip 5pt

\noindent Hence 
\begin{eqnarray*} 
V_{n+2k} & \ge & \sum_{i=1}^m {\mu}_i V_{n+2i} \\ 
& \ge & s\left( F_{2k}-F_{2m} \right) + V_{n+2m} \\
& \ge & s\left( F_{2k}-F_{2k-2} \right) + V_{n+2k-2} \\
& = & s\left( F_{2k-1} \right) + V_{n+2k-2} \\
& = & V_{n+2k} + V_n \\ 
& > & V_{n+2k}. 
\end{eqnarray*}

\noindent This contradiction completes the proof of the Proposition. 
\end{Pf}
\vskip 5pt

\begin{cor} \label{embed_dim_special}
\begin{itemize}
\item[]
\item[{\rm (i)}]
If $S=\la F_n, F_{n+2}, F_{n+4}, \ldots \ra$, $n \ge 3$, then $\e(S)=\left\lf\frac{n+1}{2}\right\rf$.
\item[{\rm (ii)}]
If $S=\la L_n, L_{n+2}, L_{n+4}, \ldots \ra$, $n \ge 4$, then $\e(S)=\left\lf\frac{n+3}{2}\right\rf$. 
\end{itemize}
\end{cor}

\begin{Pf}
This is a direct application of Proposition \ref{embed_dim}. 
\begin{itemize}
\item[{\rm (i)}]
Note that $F_{2\k-2} \le F_{2\k-1}-1<F_{2\k}$ and $F_{2\k-2}<F_{2\k}-1<F_{2\k}$. Thus, $F_n-1 \in \big[F_{2\k-2},F_{2\k}\big)$ both when $n=2\k-1$ and when $n=2\k$. Hence $\k=\left\lf\frac{n+1}{2}\right\rf$. 
\item[{\rm (ii)}]
Since $L_n=F_{n+1}+F_{n-1}$, we have $F_{2\k} \le L_{2\k-1}-1=F_{2\k}+F_{2\k-2}-1<F_{2\k+2}$ and $F_{2\k} \le L_{2\k}-1=F_{2\k+1}+F_{2\k-1}-1<F_{2\k+2}$. Thus, $L_n-1 \in \big[F_{2\k},F_{2\k+2}\big)$ both when $n=2\k-1$ and when $n=2\k$. Hence $\k=\left\lf\frac{n+3}{2}\right\rf$. 
\end{itemize}
\end{Pf}
\vskip 5pt

\noindent Computation of $\g(S)$ is difficult in the general case. In the following result we compute the genus in the special case of Fibonacci subsequences. The result is in terms of the $k^{\text{th}}$ term of  sequence that satisfies a second order recurrence, and that can be explicitly solved. 
\vskip 5pt

\begin{thm} \label{genus_special}
Let $\{{\s}_k\}$ be the second order recurrence given by 
\[ {\s}_k = 3{\s}_{k-1} - {\s}_{k-2} + F_{n+4k-1} + F_{2k-1} F_{n+2k}, \quad k \ge 2, \]
with ${\s}_0=0$ and ${\s}_1=3F_{n+2}$. Let $\{{\r}_k\}$ be the related second order recurrence given by 
\[ {\r}_k = {\s}_k - 2{\s}_{k-1} - F_{2k+1}F_{n+2k}, \quad k \ge 2.  \]
If $S=\la F_n, F_{n+2}, F_{n+4}, \ldots \ra$ and $\e(S)=\k$, then  
\[ \g(S) = \begin{cases} 
                 \ds\frac{{\s}_{\k-2}+{\r}_{\k-1}}{F_n} - \ds\frac{F_n-1}{2} & \mbox{ if $n=2\k-1$}; \\[8pt]
                 \ds\frac{{\s}_{\k-1}}{F_n} - \ds\frac{F_n-1}{2} & \mbox{ if $n=2\k$}.
                \end{cases}
\]
\end{thm}

\begin{Pf}
Applying Theorem \ref{Frob_gen}, part (ii) to $V_n=F_n$, we must show 
\[ \sum_{x=1}^{F_n-1} s(x) = \begin{cases}
                                                {\s}_{\k-2}+{\r}_{\k-1} & \mbox{ if $n=2\k-1$}; \\
                                                {\s}_{\k-1} & \mbox{ if $n=2\k$}.
                                              \end{cases}
\] 
\vskip 5pt

\noindent For positive integer $k$ and $\ell \in \{1,2\}$, define 
\begin{equation} \label{newsums}
{\r}_{k,\ell} = \sum_{x=F_{2k}}^{F_{2k+\ell}-1} s(x), \quad {\s}_k = \sum_{x=1}^k {\r}_{x,2} = \sum_{x=1}^{F_{2k+2}-1} s(x). 
\end{equation}

\noindent Fix $x \in \{1,\ldots,F_n-1\}$, and let $k$ be defined by $F_{2k} \le x < F_{2k+2}$. Then ${\l}_k \ge 1$ by the Greedy Algorithm, so that $s(x)=s(x-F_{2k})+F_{n+2k}$. Hence  
\[ {\r}_{k,2} = \sum_{x=F_{2k}}^{F_{2k+2}-1} \big( s(x-F_{2k})+F_{n+2k} \big) = \left( \sum_{x=1}^{F_{2k+1}-1} s(x) \right) + F_{2k+1} F_{n+2k} = {\s}_{k-1} + {\r}_{k,1} + F_{2k+1} F_{n+2k}, \]
so that 
\begin{equation} \label{eqn1}
{\s}_k = {\s}_{k-1} + {\r}_{k,2} = {\s}_{k-1} + {\s}_{k-1} + {\r}_{k,1} + F_{2k+1} F_{n+2k} = 2{\s}_{k-1} + {\r}_{k,1} + F_{2k+1} F_{n+2k}. 
\end{equation} 
Similarly
\begin{equation} \label{eqn2}
{\r}_{k,1} = \sum_{x=F_{2k}}^{F_{2k+1}-1} \big( s(x-F_{2k})+F_{n+2k} \big) = \left( \sum_{x=1}^{F_{2k-1}-1} s(x) \right) + F_{2k-1} F_{n+2k} = {\s}_{k-2} + {\r}_{k-1,1} + F_{2k-1} F_{n+2k}. 
\end{equation}
Replacing $k$ by $k-1$ in eqn.~\eqref{eqn1}, we get
\begin{equation} \label{eqn3}
{\s}_{k-1} = 2{\s}_{k-2} + {\r}_{k-1,1} + F_{2k-1} F_{n+2k-2}. 
\end{equation} 
From eqns.~\eqref{eqn1}, \eqref{eqn2}, \eqref{eqn3}, and using Proposition \ref{V_identity}, part (i), we have
\begin{eqnarray*} 
{\s}_k - {\s}_{k-1} & = & 2{\s}_{k-1} -  2{\s}_{k-2} + \big( {\r}_{k,1} - {\r}_{k-1,1} \big) + \big( F_{2k+1} F_{n+2k} - F_{2k-1} F_{n+2k-2} \big) \\
& = & 2{\s}_{k-1} - {\s}_{k-2} +  F_{2k-1} F_{n+2k} + \big( F_{2k+1} F_{n+2k} + F_{2k} F_{n+2k-1} \big) \\
& & - \big( F_{2k} F_{n+2k-1} + F_{2k-1} F_{n+2k-2} \big) \\
& = & 2{\s}_{k-1} - {\s}_{k-2} +  F_{2k-1} F_{n+2k} + F_{n+4k} - F_{n+4k-2} \\
& = & 2{\s}_{k-1} - {\s}_{k-2} +  F_{2k-1} F_{n+2k} + F_{n+4k-1},
\end{eqnarray*}
so that 
\begin{equation} \label{sigma}
{\s}_k = 3{\s}_{k-1} - {\s}_{k-2} +  F_{2k-1} F_{n+2k} + F_{n+4k-1}. 
\end{equation}
\vskip 5pt
 
\noindent Thus, $\{{\s}_k\}$ satisfies the recurrence given in the Theorem. Moreover, ${\s}_0=0$ and ${\s}_1=s(1)+s(2)=F_{n+2}+2F_{n+2}=3F_{n+2}$. 
\vskip 5pt

\noindent Now  
\[ {\r}_{k,1} = {\s}_k - 2{\s}_{k-1} - F_{2k+1} F_{n+2k} \]
can be used to determine ${\r}_{k,1}$ from the sequence $\{{\s}_k\}$; this is the recurrence for the sequence $\{{\r}_k\}$ with ${\r}_k$ being used for ${\r}_{k,1}$ here.  
\vskip 5pt

\noindent Recall that $\e(S)=\k$ where $F_{2(\k-1)} \le F_n-1 < F_{2\k}$ by Theorem \ref{embed_dim}. If $n=2\k$, then $\sum_{x=1}^{F_n-1} s(x) = \sum_{x=1}^{F_{2\k}-1} s(x) = {\s}_{\k-1}$. If $n=2\k-1$, then $\sum_{x=1}^{F_n-1} s(x) = \sum_{x=1}^{F_{2\k-1}-1} s(x) = {\s}_{\k-2}+{\r}_{k-1,1}$. 

\noindent This completes the proof of the Theorem. 
\end{Pf}
\vskip 5pt

\begin{rem}
Theorem \ref{genus_special} provides a formula for $\g(S)$ when $S$ is generated by the Fibonacci subsequence. This formula is in terms of ${\s}_k$ and ${\r}_k$, both of which can be determined in closed form because the recurrence has constant coefficients and the non-homogenous part is a linear combinations of powers of distinct real numbers. However, this expression is not particularly simple and so has been omitted here.   
\end{rem}
\vskip 5pt

\noindent {\bf Acknowledgement.} This work, carried out by the first two authors under the supervision of the third author over the period May -- October 2021, formed the basis for the Summer Undergraduate Research Award at the Indian Institute of Technology Delhi. The third author acknowledges the contribution of Ryan Azim Shaikh in thoroughly reading this manuscript and simplifying portions therein.

\end{document}